\documentclass[12pt]{article}

\usepackage{ amsmath, amssymb, color}
\usepackage{enumerate}

\newtheorem{theorem}{\bf Theorem}[section]

\newtheorem{lemma}[theorem]{\bf Lemma}

\newcommand{\proof}{\noindent{\bf Proof.\ }}
\newcommand{\qed}{\hfill $\square$ \bigskip}

\newcommand{\rad}{{\rm rad}}
\newcommand{\diam}{{\rm diam}}

\begin{document}

\title{On the difference between the eccentric connectivity index and eccentric distance sum of graphs}

\author{
Yaser Alizadeh $^{a}$ \and  Sandi Klav\v zar $^{b}$}

\date{}

\maketitle

\begin{center}
$^a$ Department of Mathematics, Hakim Sabzevari University, Sabzevar, Iran\\
e-mail: {\tt y.alizadeh.s@gmail.com} \\
\medskip
$^b$ Faculty of Mathematics and Physics, University of Ljubljana, Slovenia \\
e-mail: {\tt sandi.klavzar@fmf.uni-lj.si}

\end{center}

\begin{abstract}
The eccentric connectivity index of a graph $G$ is $\xi^c(G) = \sum_{v \in V(G)}\varepsilon(v)\deg(v)$, and the eccentric distance sum is $\xi^d(G) = \sum_{v \in V(G)}\varepsilon(v)D(v)$, where $\varepsilon(v)$ is the eccentricity of $v$, and $D(v)$ the sum of distances between $v$ and the other vertices. A lower and an upper bound on $\xi^d(G) - \xi^c(G)$ is given for an arbitrary graph $G$. Regular graphs with diameter at most $2$ and joins of cocktail-party graphs with complete graphs form the graphs that attain the two equalities, respectively. Sharp lower and upper bounds on $\xi^d(T) - \xi^c(T)$ are given for arbitrary trees. Sharp lower and upper bounds on $\xi^d(G)+\xi^c(G)$ for arbitrary graphs $G$ are also given, and a sharp lower bound on $\xi^d(G)$ for graphs $G$ with a given radius is proved.  
\end{abstract}


\medskip\noindent
{\bf AMS Subj.\ Class (2020)}: 05C12, 05C09, 05C92 


\section{Introduction}
\label{sec:intro}

In this paper we consider simple and connected graphs. If $G = (V(G), E(G))$ is a graph and $u,  v \in V(G)$, then the {\em distance} $d_G(u, v)$ between $u$ and $v$ is the number of edges on a shortest $u,v$-path. The eccentricity of a vertex and its total distance are distance properties of central interest in (chemical) graph theory; they are defined as follows. The {\em eccentricity} $\varepsilon_G(v)$ of a vertex $v$ is the distance between $v$ and a farthest vertex from $v$, and the {\em total distance} $D_G(v)$ of $v$ is the sum of distances between $v$ and the other vertices of $G$. Even more fundamental property of a vertex in (chemical) graph theory is its degree (or valence in chemistry), denoted by $\deg_G(v)$. (We may skip the index $G$ in the above notations when $G$ is clear.) Multiplicatively combining two out of these three basic invariants naturally leads to the {\em eccentric connectivity index} $\xi^c(G)$, the {\em eccentric distance sum} $\xi^d(G)$, and the {\em degree distance} $DD(G)$, defined as follows: 
\begin{align*}
\xi^c(G) & = \sum_{v \in V(G)}\varepsilon(v)\deg(v)\,. \\
\xi^d(G) & = \sum_{v \in V(G)}\varepsilon(v)D(v)\,. \\
DD(G) & = \sum_{v \in V(G)}\deg(v)D(v)\,.
\end{align*}
$\xi^c$ was introduced by Sharma, Goswami, and Madan~\cite{sha-1997}, $\xi^d$ by Gupta, Singh, and Madan~\cite{Gupta}, and $DD$ by Dobrynin and Kochetova~\cite{Dob1994} and by Gutman~\cite{Gut-1994}. These three topological indices are  well investigated, selected contrubutions to the eccentric connectivity index are~\cite{Hauweele-2019, Ill-2011, xu-2016}, to  the eccentric distance sum~\cite{chen-2018, Illic, xie-2019}, and to the degree distance~\cite{li-2015, Tom-2008, wang-2016}. The three invariants were also compared to other invariants, cf.~\cite{Dank2014, Das2011, das-2015, Das-2016, xu-2016b}. For information on additional topological indices based on eccentricity see~\cite{Madan-2010}.

In~\cite{hua-2018} the eccentric distance sum and the degree distance are compared, while in~\cite{xu-2017} the difference between the eccentric connectivity index and the (not defined here) connective eccentricity index is studied. The primary motivation for the present paper, however, are the papers~\cite{hua-2019, zhang-2019} in which $\xi^d(G) - \xi^c(G)$ was investigated. In~\cite{zhang-2019}, Zhang, Li, and Xu, besides other results on the two indices, determined sharp upper and lower bounds on $\xi^d(G) - \xi^c(G)$ for graphs $G$ of given order and  diameter $2$. Parallel results were also derived for sub-classes of diameter $2$ graphs with specified one of the minimum degree, the connectivity, the edge-connectivity, and the independence number. Hua, Wang, and Wang~\cite{hua-2019} extended the last result to general graphs. More precisely, they characterized the graphs that attain the minimum value of $\xi^d(G) - \xi^c(G)$ among all connected graphs $G$ of given independence number. They also proved a related result for connected graphs with given matching number. 

In this paper we continue the investigation along the lines of~\cite{hua-2019, zhang-2019} and proceed as follows. In the rest of this section definitions and some observations needed are listed. In Section~\ref{sec:difference-general}, we give a lower and an upper bound on $\xi^d(G) - \xi^c(G)$ and in both cases characterize the equality case. The upper bound involves the Wiener index, the first Zagreb index, as well as the degree distance of $G$. In Section~\ref{sec:trees} we focus on trees and first prove that among all trees $T$ with given order and diameter, $\xi^d(T)-\xi^c(T)$ is minimized on caterpillars. Using this result we give a lower bound on $\xi^d(T)-\xi^c(T)$ for all trees $T$ with given order, the bound being sharp precisely on stars. We also give a sharp upper bound on $\xi^d(T)-\xi^c(T)$ for trees $T$ with given order. In the last section we give a sharp lower bound and a sharp upper bound on $\xi^d(G)+\xi^c(G)$, compare $\xi^d(G)$ with $\xi^c(G)$ for graphs $G$ with not too large maximum degree, and give a sharp lower bound on $\xi^d(G)$ for graphs $G$ with a given radius.  

\subsection{Preliminaries}

The order and the size of a graph $G$ will be denoted by $n(G)$ and $m(G)$, respectively. The star of order $n\ge 2$ is denoted by $S_n$; in other words, $S_n = K_{1,n-1}$.  If $n\ge 2$, then the {\em cocktail party graph} $CP_{2n}$ is the graph obtained from $K_{2n}$ by removing a perfect matching. The {\em join} $G\oplus H$ of graphs $G$ and $H$ is the graph obtained from the disjoint union of $G$ and $H$ by connecting by an edge every vertex of $G$ with every vertex of $H$.  The maximum degree of a vertex of $G$ is denoted by $\Delta(G)$. A graph $G$ is {\em regular} if all vertices have the same degree. The {\em first Zagreb index}~\cite{gutman-1972} $M_1(G)$ of $G$ is the sum of the squares of the degrees of the vertices of $G$. The {\em Wiener index}~\cite{wiener} $W(G)$ of $G$ is the sum of distances between all pairs of vertices in $G$.  

The {\em diameter} $\diam(G)$ and the {\em radius} $\rad(G)$ of a graph $G$ are the maximum and the minimum vertex eccentricity in $G$, respectively. A graph $G$ is {\em self-centered} if all vertices have the same eccentricity. It this eccentricity is $d$, we further say that $G$ is {\em $d$-self-centered}. The {\em eccentricity} $\varepsilon(G)$ of $G$ is $$\varepsilon(G)=\sum_{v \in V(G)}\varepsilon(v)\,.$$
The eccentric connectivity index of $G$ can be equivalently written as 
\begin{equation}
\label{eq:eci}
\xi^c(G) = \sum_{uv \in E(G)}\varepsilon(u) + \varepsilon(v)\,,
\end{equation}
and the eccentric distance sum as
\begin{equation}
\label{eq:eds}
\xi^d(G) = \sum_{ \{u,v\}\subseteq V(G)}(\varepsilon(u)+\varepsilon(v))d(u,v)\,.
\end{equation}

\section{The difference on general graphs}
\label{sec:difference-general}

In this section we give some sharp upper and lower bounds on $\xi^d(G) - \xi^c(G)$ for an arbitrary graph $G$. The bounds are in terms of the eccentricity, the Wiener index, the first Zagreb index, the degree distance, the maximum degree, the size, and the order of $G$. 

\begin{theorem}
\label{thm:difference}
If $G$ is a connected graph, then the following hold.
\begin{enumerate}[(i)]
\item $\xi^d(G) - \xi^c(G) \ge 2\big(n(G)-1-\Delta(G)\big)\varepsilon(G)$. Moreover, the equality holds if and only if $G$ is a regular graph with $\diam(G) \le 2$. 
\item $\xi^d(G) - \xi^c(G)  \le 2n(G)\big(W(G)-m(G)\big) + M_1(G) - DD(G)$. Moreover, the equality holds if and only if $G\in \{P_4\} \cup \{CP_{2k}\oplus K_{n(G)-2k}:\ 0\le k\le n/2\}$.  
\end{enumerate}
\end{theorem}

\proof
(i) Let $v$ be a vertex of $G$. If $w$ is not adjacent to $v$, then $d(v,w)\ge 2$ and consequently $D(v) - \deg(v) \ge 2(n(G)-1-\Delta(G))$. Thus:
\begin{eqnarray*}
 \xi^d(G)-\xi^c(G) &=& \sum_{v \in V(G)}\varepsilon(v)\big(D(G)-\deg(v)\big)\\
 & \ge & \sum_{v \in V(G)} 2\, \varepsilon(v)\big(n(G) -1-\Delta(G)\big) \\
 & = & 2\,\varepsilon(G)\big( n(G)-1-\Delta(G)\big)\,.
 \end{eqnarray*}
The equality holds if and only if $D(v) - \deg(v) = 2(n(G)-1-\Delta(G))$ for every vertex $v$. As the last equality in particular holds for a vertex of maximum degree, we infer that $G$ must be regular. Then the condition  $D(v) - \deg(v) = 2(n(G)-1-\Delta(G))$ simplifies to  
\begin{equation}
\label{eq:D-delta}
D(v) + \Delta(G) = 2n(G)-2\,.
\end{equation} 
Suppose that $\diam(G) = d$, and let $x_i$, $i\in\{2,\ldots, d\}$, be the number of vertices at distance $i$ from $v$. Then $n(G) = 1 + \Delta(G) + x_2 + \cdots + x_d$ and $D(v) = \Delta(G) + 2x_2 + \cdots + dx_d$. Plugging these equalities into~\eqref{eq:D-delta} yields
$$2\Delta(G) + 2x_2 + \cdots + dx_d = 2 + 2\Delta(G) + 2x_2 + \cdots + 2x_d - 2$$
which implies that $x_3 = \cdots = x_d = 0$, that is, $\diam(G) = 2$. Finally, if $\diam(G) = 2$, then $D(v) = \Delta(G) + 2(n(G)-\Delta(G)-1)$, so~\eqref{eq:D-delta} is fulfilled for every regular graph of diameter $2$. Clearly,~\eqref{eq:D-delta} is also fulfilled for graphs of diameter $1$, that is, complete graphs.   

(ii) If $v \in V(G)$, then clearly $\varepsilon(v) \le n(G)-\deg(v)$. Then we deduce that 
\begin{eqnarray*}
\xi^d(G) - \xi^c(G) &=& \sum_{v \in V(G)}\varepsilon(v)\big( D(v)-\deg(v) \big) \\
& \le &  \sum_{v \in V(G)} \big( n(G)-\deg(v) \big)\big( D(v)-\deg(v) \big) \\
&=& n(G)\sum_{v \in V(G)}\big(D(v)-\deg(v) \big) + \sum_{v \in V(G)}\deg(v)^2 \\  
& & - \sum_{v \in V(G)}\deg(v)D(v) \\
&=& 2n(G)\big(W(G)-m(G)\big) + M_1(G) - DD(G)\,. 
\end{eqnarray*}
The equality in the above computation holds if and only if $\varepsilon(v) = n(G) - \deg(v)$ holds for all $v \in V (G)$. So suppose that $G$ is a graph for which $\varepsilon(v) = n(G) - \deg(v)$ holds for all $v \in V (G)$ and distinguish the following two cases. 

Suppose first that $\diam(G)\ge 3$. Let $P$ be a diametral path in $G$ and let $v$ and $v'$ be its endpoints. Since   $\varepsilon(v) = n(G) - \deg(v)$ and $|V(P) \setminus N[v]| = \varepsilon(v) - 1$, it follows that $n(G) = 1 + \deg(v) + |V(P) \setminus N[v]|$. The latter means that $V(G) = N[v]\cup V(P)$. Since $\diam(G) = \varepsilon(v) \ge 3$ it follows that $\deg(v') = 1$. Since we have also assumed that $\varepsilon(v') = n(G) - \deg(v')$ holds we see that  $\varepsilon(v') = n(G) -1$ which in turn implies that $G$ is a path. Among the paths $P_n$, $n\ge 4$, the path $P_4$ is the unique one which fulfills the condition $\varepsilon(v) = n - \deg(v)$ for all $v \in V (P_n)$.

Suppose second that $\diam(G)\le 2$. Then $\varepsilon(v)\in \{1,2\}$ for every $v\in (G)$. Since $\varepsilon(v) = n(G) - \deg(v)$ it follows that $\deg(v)\in \{n(G)-1, n(G)-2\}$. Let $V_1 = \{v:\ \deg(v) = n(G) - 1\}$ and $V_2 = \{v:\ \deg(v) = n(G) - 2\}$. Then $V(G) = V_1 \cup V_2$. Clearly, the subgraph of $G$ induced by $V_1$ is complete, and there are all possible edges between $V_1$ and $V_2$. Moreover, the complement of the subgraph of $G$ induced by $V_2$ is a disjoint union of copies of $K_2$, which means that $V_2$ induces a cocktail party graph. In summary, $G$ must be of the form $CP_{2k}\oplus K_{n(G)-2k}$, where $0\le k\le n/2$. On the other hand, the condition $\varepsilon(v) = n(G) - \deg(v)$ clearly holds for each vertex of $CP_{2k}\oplus K_{n(G)-2k}$, hence these graphs together with $P_4$ from the previous case are precisely the graphs that attain the equality. 
\qed

\section{The difference on trees}
\label{sec:trees}

In this section we turn our attention to  $\xi^d(T) - \xi^c(T)$ for trees $T$, and in particular on extremal trees regarding this difference. 

\begin{theorem} 
\label{caterpilar}
Among all trees $T$ with given order and diameter, $\min\{ \xi^d(T)-\xi^c(T)\}$ is achieved on caterpillars.
\end{theorem}

\proof
Fix the order and diameter of trees to be considered. Let $T$ be an arbitrary tree that is not a caterpillar with this fixed order and diameter. Let $P$ be a diametral path of $T$ connecting $x$ to $y$. Then the eccentricity of each vertex $w$ of $T$ is equal to $\max\{d(w,x), d(w,y)\}$. Let $z\ne x,y$ be a vertex of $P$ and let $T_z$ be a maximal subtree of $T$ which contains $z$ but no other vertex of $P$. We may assume that $z$ can be selected such that $\varepsilon_{T_z}(z) = k \ge 2$, for otherwise $T$ is a caterpillar. Let $u$ be vertex of $T_z$ with $d(u,z)=k-1$ and let $v$ be the neighbor of $u$ with $d(v,z)=k-2$. Let $S = N(u)\setminus \{v\}$ and let $s = |S|$. Note that $s > 0$. Let now $T'$ be the tree obtained from $T$ by replacing the edges between $u$ and the vertices of $S$ with the edges between $v$ and the vertices of $S$. 

\medskip\noindent
{\bf Claim A}: $\xi^d(T) - \xi^c(T) > \xi^d(T') - \xi^c(T')$.\\
Set $X_d = \xi^d(T) - \xi^d(T')$ and $X_c = \xi^c(T) - \xi^c(T')$. To prove the claim it is equivalent to show that  $X_d - X_c > 0$. 

For a vertex $w \in V(G)\setminus (S \cup \{u\})$ we have $D_{T'}(w) = D_{T}(w) - s$ and $\varepsilon_{T'}(w)\le \varepsilon_{T}(w)$. Moreover if $w \in S$, then $\varepsilon_{T'}(w)= \varepsilon_{T}(w) -1$ and $D_T(w) = D_{T'}(w) + n-s-2$. With these facts in hand we can compute as follows. 
\begin{eqnarray*}
X_d &=& \sum_{w \in V(T)}\varepsilon_T(w)D_T(w) -  \sum_{w \in V(T')}\varepsilon_{T'}(w)D_{T'}(w) \\
&=& \varepsilon_{T}(u)D_{T}(u)- \varepsilon_{T'}(u)D_{T'}(u) + \varepsilon_{T}(v)D_{T}(v)-\varepsilon_{T'}(v)D_{T'}(v)\\
& &+  \sum_{w \in S}\varepsilon_T(w)D_T(w) - \varepsilon_{T'}(w)D_{T'}(w) \\ &&+ \sum_{w \in V(T)- ( S\cup \{u,v\})}\varepsilon_T(w)D_T(w) - \varepsilon_{T'}(w)D_{T'}(w)\\
&\ge& s(\varepsilon_T(v)-\varepsilon_T(u)) 
+ \sum_{w \in V(T)- ( S\cup \{u,v\})} \varepsilon_T(w)s \\&& +  \sum_{w \in S}\big( \varepsilon_T(w)D_T(w) - (\varepsilon_T(w)-1)(D_T(w)-n+2+s) \big) \\
&=& -s + \sum_{w \in V(T)- ( S\cup \{u,v\})} \varepsilon_T(w)s   \\ && +  \sum_{w \in S}\big( (D_T(w)-n+2+s ) - \varepsilon_T(w)(-n+2+s) \big) \\
&=& -s + \sum_{w \in V(T)\setminus ( S\cup \{u,v\})} \varepsilon_T(w)s + (n-s-2)\sum_{w \in S}\varepsilon_T(w)-1 + D_T(w)\\
&=&  -s + \sum_{w \in V(T)\setminus ( S\cup \{u,v\})} \varepsilon_T(w)s \\ & & + s(n-s-2)\varepsilon_T(u) + s(D_T(u)+ n-2)\\
&=& s\big[ \varepsilon(T)-\varepsilon_T(u)(s+2)-s+1   +  (n-s-2)\varepsilon_T(u)+ D_T(u)+ n-3) \big]\,.
\end{eqnarray*}
Similarly, but shorter, we get that $X_c = 2s$. Thus 
\begin{align*}
X_d - X_c  & \ge s\big[ \varepsilon(T)-\varepsilon_T(u)(s+2) \\ 
 &\quad +  (n-s-2)\varepsilon_T(u)+ D_T(u)+ n-s-4) \big] \\
 & > 0\,.
\end{align*}
This proves Claim A. If $T'$ is not a caterpillar, we can repeat the construction as many times as required to arrive at a caterpillar. Since at each step the value of $\xi^d - \xi^c$ is decreased, the minimum of this difference is indeed achieved on caterpillars. 
\qed

\begin{theorem}\label{min}
If $T$ is a tree of order $n\ge 3$, then
\[\xi^d(T) - \xi^c(T)\ge  4 n^2-12 n+8\,.\]
Moreover, equality holds if and only if $T = S_n$.
\end{theorem}

\proof
Let $n\ge 3$ be a fixed integer. By Theorem~\ref{caterpilar}, it suffices to consider caterpillars. More precisely, let $T$ be a caterpillar of order $n$ and with $\diam(T) = d \ge 3$. Then we wish to prove that $\xi^d(T) - \xi^c(T) > \xi^d(S_n) - \xi^c(S_n) =  4 n^2-12 n+8$. The latter equality is straightforward to check, for the strict inequality we proceed as follows. 

Let $w, z\in V(T)$ be two adjacent vertices of eccentricities $d-1$ and $d-2$, respectively. Let $S = N(w)\setminus\{z\}$ and set $s = |S|$. As $\varepsilon(w) = d-1$, we have $s\ge 1$. Let further $S_1 = V(G) \setminus (S\cup \{w,z\})$. Construct now a tree $T'$ from $T$ by replacing the edges between $w$ and the vertices of $S$ with the edges between $z$ and the vertices of $S$.  Note that $\deg_T(w) = \deg_{T'}(w)+s = 1+s$ and $\deg_T(z) = \deg_{T'}(z)- s$, while the other vertices have the same degree in $T$ and $T'$. Further, it is straightforward to verify the following relations:
\begin{align*}
D_T(w) & = D_{T'}(w)- s,\quad \varepsilon_T(w)=  \varepsilon_{T'}(w); \\
D_T(z) & = D_{T'}(z)+ s,\quad \varepsilon_{T'}(z)\le\varepsilon_T(z)\le  \varepsilon_{T'}(z)+1; \\
D_T(x) & = D_{T'}(x) + n-s-2,\quad \varepsilon_T(x)=\varepsilon_{T'}(x)+1\ (x \in S); \\
D_T(y) & = D_{T'}(y) + s,\quad \varepsilon_{T'}(y)\le \varepsilon_T(y)\le \varepsilon_{T'}(y) + 1\ (y \in S_1).\\
\end{align*}
Setting $X_d = \xi^d(T)-\xi^d(T')$ we have: 
\begin{eqnarray*}
X_d & = & \sum_{v \in \{w,z\}}D_T(v)\varepsilon_{T}(v)-D_{T'}(v)\varepsilon_{T'}(v) + \sum_{v \in S}D_{T}(v)\varepsilon_{T}(v)-D_{T'}(v)\varepsilon_{T'}(v)\\&&+ \sum_{v \in S_1 }D_T(v)\varepsilon_{T}(v)-D_{T'}(v)\varepsilon_{T'}(v) \\
&\ge & s(\varepsilon_{T'}(z)-\varepsilon_T(w)) + \sum_{v \in S}D_{T}(v)\varepsilon_{T}(v)- (D_{T}(v) -(n-s-2))(\varepsilon_{T}(v)-1)\\
&&+ \sum_{v \in S_1 }D_T(v)\varepsilon_{T}(v)-(D_{T}(v)-s)\varepsilon_{T}(v) \\
&\ge& -s + (n-s-2)\sum_{v \in S}\varepsilon_T(v)+ \sum_{v \in S}D_T(v) - s(n-s-2) + s \sum_{v \in S_1}\varepsilon_T(v)\\
&\ge & -s+ 3s(n-s-2) + s(2(n-s-2)+2s+1) - s(n-s-2) \\ 
& & + 3s(n-s-2)\\
&=& 5s(n-s-3) + 2s(n-1)\,.
\end{eqnarray*}
Similarly, setting $X_c = \xi^c(T)-\xi^c(T')$, we have
\begin{eqnarray*}
X_c &=& \sum_{v \in \{w,z\}} \big(\deg_T(v)\varepsilon_{T}(v)-\deg_{T'}(v)\varepsilon_{T'}(v)\big)\\
&&+\sum_{v \in S} \big( \deg_{T}(v)\varepsilon_{T}(v)-\deg_{T'}(v)\varepsilon_{T'}(v)\big)\\
&&+ \sum_{v \in S_1 } \big(\deg_T(v)\varepsilon_{T}(v)-\deg_{T'}(v)\varepsilon_{T'}(v)\big) \\
&\le & s\varepsilon_T(w) + \deg_T(z)\varepsilon_T(z) - (\deg_T(z)+s)(\varepsilon_T(z)-1) \\
&& + s + 
\sum_{v \in S_1 }\deg_T(v)\varepsilon_{T}(v)-\deg_{T}(v)(\varepsilon_{T}(v)-1)\\
&=& 2s + \deg_T(z) + \sum_{v \in S_1}\deg(v)\\ &=& 2n -3.
\end{eqnarray*}
Therefore, 
\begin{eqnarray*}
X_d - X_c \ge \big(5s(n-s-3) + 2s(n-1)\big) -\big(2n -3\big) > 0\,,
\end{eqnarray*}
that is, $\xi^d(T)-\xi^c(T) > \xi^d(T')-\xi^c(T')$. Repeating the above transformation until $S_n$ is constructed finishes the argument. 
\qed

To bound the difference $\xi^d(T) - \xi^c(T)$ for an arbitrary tree $T$ from above, we first recall the following result. 

\begin{lemma}\label{illic} \cite[Theorem 2.1]{Illic}
 Let $w$ be a vertex of graph $G$. For non-negative integers $p$ and $q$, let $G(p,q)$ denotes the graph obtained from $G$ by attaching to vertex $w$ pendant paths $P = wv_1\cdots v_p$ and $Q = wu_1 \cdots u_q$ of lengths $p$ and $q$, respectively. Let $G(p+q,0)=G(p,q)- wu_1+v_pu_1$. If $r = \varepsilon_G(w)$ and $r\ge p\ge q \ge 1$, then 
 \begin{align*}
 \xi^d(G(p+q,0))- \xi^d(G(p,q)) \ge &\frac{pq}{6}\big[ 6D_G(w)+p(2p-3) + q(2q-3) + 3pq -12r  \\&+ 6n(G)(p+q+r+1) + 6\sum_{v \in V(G)}\varepsilon(v)\big]\,. 
 \end{align*}
\end{lemma}

\begin{lemma}\label{xic}
Let $G$,  $p$, $q$, $G(p,q)$, and $G(p+q,0)$ be as in Lemma~\ref{illic}. Then
$$\xi^c(G(p+q,0)) - \xi^c(G(p,q)) \le q(3p + 2m(G)-1)\,.$$
\end{lemma}

\proof 
Let $\deg'(v)$ and $\varepsilon'(v)$ (resp.\ $\deg(v)$ and $\varepsilon(v)$) denote the degree and the eccentricity of $v$ in $G(p+q,0)$ (resp. $G(p,q)$). Then we have:  
\begin{align*}
\deg'(w) &= \deg(w)-1, \quad \varepsilon'(w) \le \varepsilon(w) + q; \\
\deg'(v_i)&=\deg(v_i),  i\in [p-1],\quad  \deg'(v_p)=\deg(v_p)+1;   \\
\varepsilon'(v_i) &\le \varepsilon(v_i) + q, i\in [p]; \\
\deg'(u_j)&=\deg(u_j), \quad \varepsilon'(u_j) = \varepsilon(u_j)+ p; \\
\varepsilon'(x)&\le \varepsilon(x) + q,  x \in V(G)\,.
\end{align*}
Moreover, the degrees of vertices in $G(p+q,0)$ do not decrease. Calculating the difference of contributions of vertices in $\xi^c$ for $G(p+q,0)$ and $G(p,q)$, we can estimate the difference $X_c =  \xi^c(G(p+q,0)) - \xi^c(G(p,q))$ as follows: 
\begin{eqnarray*}
X_c & \le & \sum_{w\neq x \in V(G)}\deg(x) q  + \sum_{i=1}^{q} \deg(u_i)p + \sum_{i=1}^{ \lfloor \frac{p}{2} \rfloor}\deg(v_i)q\\
& & +  \varepsilon(v_p) + (\deg(w)-1)(r+q)-\deg(w)r\\
& = & \big(2m(G)-\deg(w) \big)q + (2q-1)p + pq +p +q(\deg(w)-1) \\
& = & 2qm(G) + 3pq - q\,. 
\end{eqnarray*}
\qed
    
\begin{theorem}
If $T$ is a tree of order $n$, then  
$$\xi^d(T) - \xi^c(T) \le \begin{cases} \frac{25 n^4}{96}-\frac{n^3}{6}-\frac{89 n^2}{48}+\frac{19 n}{6}-\frac{45}{32}; & n\ odd\,, \\[5pt]  \frac{25 n^4}{96}-\frac{n^3}{6}-\frac{43 n^2}{24}+\frac{19 n}{6}-2;  &  n\ even\,.
\end{cases}$$ 
Moreover, equality holds if and only if $T=P_n$.
\end{theorem}

\proof
The right side of the above inequality is equal to $\xi^d(P_n)-\xi^c(P_n)$. (The value of $\xi^d(P_n)$ has been determined in~\cite{Illic}, while it is straightforward to deduce $\xi^c(P_n)$. Combining the two formulas, the polynomials from the right hand side of the inequality are obtained.) Suppose now that $T \neq P_n$. Then there is always a vertex $w$ of degree at least $3$ such that we can apply Lemmas~\ref{illic} and~\ref{xic}. Setting
$$X_{dc} = \big(\xi^d(T(p+q,0) - \xi^c(p+q,0)\big)-\big(\xi^d(T(p,q) - \xi^c(p,q) \big)$$
we have: 
\begin{eqnarray*}
X_{dc} & \ge& pqD_T(w) + \frac{pq}{6} \big( p(2p-3) + q(2q-3) \big)+  \frac{(pq)^2}{2}  -2pqr  \\
& & + pqn(T)(p+q+r+1) +pq\sum_{v \in V(T)}\varepsilon(v) - \big( 2qm(T) + 3pq - q \big)\\
& = & pq\big(D_T(w) + \sum_{v \in V(T)}\varepsilon(v) -3 \big) + \frac{pq}{6}\big(2p^2-3p+2q^2-3q+3pq \big)\\
& & + pqr(n(G)-2) + q\big(pn(T)(p+q+r)-2m(T)+1\big) > 0
\end{eqnarray*}
and the result follows. 
\qed

\section{Further comparison}
\label{sec:further}

In this concluding section we give sharp lower and upper bounds on $\xi^d(G)+\xi^c(G)$, compare $\xi^d(G)$ with $\xi^c(G)$ for graphs $G$ with $\Delta(G) \le \frac{2}{3}(n-1)$, and give a sharp lower bound on $\xi^d(G)$ for graphs $G$ with a given radius.

\begin{theorem}
If $G$ is a connected graph, then the following hold.
\begin{enumerate}[(i)]
\item $\xi^d(G) + \xi^c(G)\le 2(n(G)-1) \varepsilon(G) + 2 \diam(G) \big( W(G) + m(G) - 2{n(G) \choose 2}   \big)$.
\item $\xi^d(G) + \xi^c(G)\ge 2(n(G)-1)\varepsilon(G) + 2 \rad(G) \big( W(G) + m(G) - 2{n(G) \choose 2}   \big) $.
\end{enumerate}
Moreover, each of the equalities holds if and only if $G$ is a self-centered graph.
\end{theorem}

\proof (i) Partition the pairs of vertices of $G$ into neighbors and non-neighbors, and using~\eqref{eq:eci}, we can compute as follows:   
\begin{eqnarray*}
\xi^d(G) &=& \sum_{ \{u,v\}\subseteq V(G)}(\varepsilon(u)+\varepsilon(v))d(u,v) \\
&=& \sum_{ uv\in E(G)}(\varepsilon(u)+\varepsilon(v))  +2  \sum_{ \{u,v\}\subseteq V(G) \atop d(u,v)\ge 2} (\varepsilon(u)+\varepsilon(v))  \\ & & +    \sum_{ \{u,v\}\subseteq V(G) \atop  d(u,v)\ge 2} \big(\varepsilon(u)+\varepsilon(v)\big)\big(d(u,v)-2\big) \\ 
&=& \xi^c(G) + \sum_{\{u,v\}\subseteq V(G)}(\varepsilon(u) + \varepsilon(v)) - 2 \xi^c(G) \\
&& +\sum_{\{u,v\} \subseteq V(G)\atop d(u,v)\ge 2} \big(\varepsilon(u)+\varepsilon(v)\big)\big(d(u,v)-2\big) \\
& \le & -\xi^c(G) + 2(n(G)-1)\varepsilon(G) \\
&& + 2\diam(G) \big( W(G) + m(G) - 2 \tbinom{n(G)}{2} \big)\,.
\end{eqnarray*}
The inequality above becomes equality if and only if $\varepsilon(v) = \diam(G)$ for every $v \in V(G)$. That is, the equality holds if and only if $G$ is a self-centered graph. 

(ii) This inequality as well as its equality case are proved along the same lines as (i). The only difference is that the inequality $\varepsilon(u)+\varepsilon(v) \le 2\diam(G)$ is replaced by $\varepsilon(u)+\varepsilon(v) \ge 2\rad(G)$. 
\qed

In our next result we give a relation between $\xi^d(G)$ and $\xi^c(G)$ for graph $G$ with maximum degree at most $\frac{2}{3}(n(G)-1)$. 

\begin{theorem}
If $G$ is a graph with $\Delta(G) \le \frac{2}{3}(n-1)$, then $\xi^d(G) \ge 2\xi^c(G)$. Moreover, the equality holds if and only if $G$ is $2$-self-centered, $\frac{2}{3}(n(G)-1)$-regular graph.
\end{theorem}

\proof
Set $n = n(G)$ and let $v$ be a vertex of $G$. Since $\deg(v) < n-1$ we have  $\varepsilon(v) \ge 2$. Therefore $D(v) \ge 2(n-1) - \deg(v) $ with equality holding if and only if $\varepsilon(v) = 2$. Using the assumption that $\deg(v) \le \frac{2}{3}(n-1)$, equivalently, $ 2n-2 \ge 3\deg(v)$, we infer that $ \varepsilon(v)D(v) \ge 2\varepsilon(v)\deg(v)$. Summing over all vertices of $G$ the inequality is proved. Its derivation also reveals that the equality holds if and only if $\deg(v) =\frac{2}{3}(n-1)$ and $\varepsilon(v)=2$ for each vertex $v \in V(G)$. 
\qed

To conclude the paper we give a lower bound on the eccentric distance sum in terms of the radius of a given graph. Interestingly, the cocktail-party graphs are again among the extreme graphs. 

\begin{theorem}
If $G$ is a graph with $\rad(G) = r$, then
$$\xi^d(G)\ge \big(n(G)-1 + \tbinom{r}{2}\big)\varepsilon(G)\,.$$
Equality holds if and only if $G$ is a complete graph or a cocktail-party graph. 
\end{theorem}

\proof
Set $n = n(G)$ and let $v\in V(G)$. Let $P$ be a longest path starting in $v$. Separately considering the neighbors of $v$, the last $\varepsilon(v) - 2$ vertices of $P$, and all the other vertices, we can estimate that 
 \begin{eqnarray*}
 D(v)&\ge& \deg(v) + (3 + \cdots + \varepsilon(v)) + 2\big(n - 1 - \deg(v) - (\varepsilon(v) -2) \big)  \\
 &=& 2n - \deg(v) + \frac{\varepsilon(v)^2-3\varepsilon(v)}{2}-1\,.  
 \end{eqnarray*}
Since $n-\deg(v) \ge \varepsilon(v)$ for every vertex $v \in V(G)$, we have $D(v) \ge n+ \varepsilon(v) + \frac{\varepsilon(v)^2 - 3\varepsilon(v)}{2}-1$. Consequently, having the fact $\varepsilon(v)\ge r$ in mind, we get $D(v)\ge n-1 + {r \choose 2}$. Multiplying this inequality by $\varepsilon(v)$ and summing over all vertices of $G$ the claimed inequality is proved. 

From the above derivation we see that the equality can holds only if $\varepsilon(v)=r=n-\deg(v)$ holds for every $v\in V(G)$.  From the equality part of the proof of Theorem~\ref{thm:difference}(ii) we know that this implies $\diam(G)\le 2$. For the equality we must also have  $D(v) = n-1 + {r \choose 2}$ for every $v$. If $r = 2$ this means that $D(v) = n$ and hence $\deg(v) = n-2$. It follows that $G$ is a cocktail-party graph. And if $r = 2$, then we get a complete graph. 
\qed

\section*{Acknowledgements}

Sandi Klav\v{z}ar acknowledges the financial support from the Slovenian Research Agency (research core funding P1-0297 and projects J1-9109, J1-1693, N1-0095).

\end{document}